\documentclass{amsart}%
\usepackage{amscd}
\usepackage{amsmath}
\usepackage{amsfonts}
\usepackage{amssymb}
\usepackage{graphicx}%
\providecommand{\U}[1]{\protect\rule{.1in}{.1in}}
\newtheorem{theorem}{Theorem}
\theoremstyle{plain}
\newtheorem{acknowledgement}{Acknowledgement}

\newtheorem{definition}{Definition}

\newtheorem{lemma}{Lemma}

\newtheorem{proposition}{Proposition}

\numberwithin{equation}{section}

\begin{document}
\title[Maximal operators Lorentz]{Sharp Lorentz estimates for dyadic-like maximal operators and related Bellman functions}
\author{Antonios D. Melas}
\author{Eleftherios N. Nikolidakis}
\address{Department of Mathematics, University of Athens, Panepistimiopolis 15784,
Athens, Greece}
\email{amelas@math.uoa.gr, lefteris@math.uoc.gr}
\date{March 25, 2014}

\begin{abstract}
We precisely evaluate Bellman type functions for the dyadic maximal opeator
$\mathbb{R}^{n}$ and of maximal operators on martingales related to local
Lorentz type estimates. Using a type of symmetrization principle, introduced
for the dyadic maximal operator in earlier works of the authors we precisely
evaluate the supremum of the Lorentz quasinorm of the maximal operator on a
function $\phi$ when the integral of $\phi$ is fixed and also the same Lorentz
quasinorm of $\phi$ is fixed. Also we find the corresponding supremum when the
integral of $\phi$ is fixed and several weak type conditions are given.

\begin{acknowledgement}
This research has been co-financed by the European Union and Greek national
funds through the Operational Program "Education and Lifelong Learning" of the
National Strategic Reference Framework (NSRF). ARISTEIA I, MAXBELLMAN 2760,
research number 70/3/11913.

\end{acknowledgement}

\begin{acknowledgement}
The authors would like to thank Prof. D. Cheliotis for his help.

\end{acknowledgement}
\end{abstract}
\maketitle

\section{Introduction}

The dyadic maximal operator on $\mathbb{R}^{n}$ is defined by
\begin{equation}
M{}\,_{d}\phi{}(x)=\sup\left\{  \frac{1}{\left\vert Q\right\vert }\int
_{Q}\left\vert \phi(u)\right\vert du:x\in Q\text{, }Q\subseteq\mathbb{R}%
^{n}\text{ is a dyadic cube}\right\}  \label{i1}%
\end{equation}
for every $\phi\in L_{\text{loc}}^{1}(\mathbb{R}^{n})$ where the dyadic cubes
are the cubes formed by the grids $2^{-N}\mathbb{Z}^{n}$ for $N=0,1,2,...$.

As it is well known it satisfies the following weak type $(1,1)$ inequality
\begin{equation}
\left\vert \{x\in\mathbb{R}^{n}:M\,_{d}\phi(x)>\lambda\}\right\vert \leq
\frac{1}{\lambda}\int_{\{M\,_{d}\phi>\lambda\}}\left\vert \phi(u)\right\vert
du\text{.} \label{i2}%
\end{equation}

\noindent for every $\phi\in L^{1}(\mathbb{R}^{n})$ and every $\lambda>0$ from
which it is easy to get the following $L^{p}$ inequality
\begin{equation}
\left\Vert M_{d}\phi\right\Vert _{p}\leq\dfrac{p}{p-1}\left\Vert
\phi\right\Vert _{p} \label{i3}%
\end{equation}
for every $p>1$ and every $\phi\in L^{p}(\mathbb{R}^{n})$ which is best
possible (see \cite{Burk1}, \cite{Burk2} for the general martingales and
\cite{Wang} for dyadic ones).

An approach for studying such maximal operators is the introduction of the so
called Bellman functions (see \cite{Naz}) related to them which reflect
certain deeper properties of them by localizing. Such functions related to the
$L^{p}$ inequality (\ref{i3}) have been precisely evaluated in\ \cite{Mel1}.
Actually defining for any $p>1$
\begin{equation}
\mathcal{B}_{p}(F,f,L)=\sup\left\{  \dfrac{1}{\left\vert Q\right\vert }%
\int_{Q}(M_{d}\phi)^{p}:\operatorname{Av}_{Q}(\phi^{p})=F,\operatorname{Av}%
_{Q}(\phi)=f,\sup\limits_{R:Q\subseteq R}\operatorname{Av}_{R}(\phi
)=L\right\}  \label{i4}%
\end{equation}
where $Q$ is a fixed dyadic cube, $R$ runs over all dyadic cubes containing
$Q$, $\phi$ is nonnegative in $L^{p}(Q)$ and the variables $F,f,L$ satisfy
$0\leq f\leq L,f^{p}\leq F$ which is independent of the choice of $Q$ (so we
may take $Q=[0,1]^{n}$) it has been shown in \cite{Mel1} that
\begin{equation}
\mathcal{B}_{p}(F,f,L)=\left\{
\begin{array}
[c]{cc}%
F\omega_{p}\left(  \tfrac{pL^{p-1}f-(p-1)L^{p}}{F}\right)  ^{p} & \text{if
}L<\tfrac{p}{p-1}f\\
L^{p}+(\tfrac{p}{p-1})^{p}(F-f^{p})\text{ } & \text{if }L\geq\tfrac{p}%
{p-1}f\text{.}%
\end{array}
\right.  \label{i5}%
\end{equation}
where $\omega_{p}:$ $[0,1]\rightarrow\lbrack1,\frac{p}{p-1}]$ is the inverse
function of $H_{p}(z)=-(p-1)z^{p}+pz^{p-1}$. Actually (see \cite{Mel1}) the
more general approach of defining Bellman functions with respect to the
maximal operator on a nonatomic probability space $(X,\mu)$ equipped with a
tree $\mathcal{T}$ (see Section 2) can be taken and the corresponding Bellman
function is always the same.

There are several other problems in Harmonic Analysis where Bellman functions
naturally arise. Such problems (including the dyadic Carleson imbedding and
weighted inequalities) are described in \cite{Naz2} (see also \cite{Naz},
\cite{Naz1}) and also connections to Stochastic Optimal Control are provided,
from which it follows that the corresponding Bellman functions satisfy certain
nonlinear second order PDE.

The exact computation of a Bellman function is a difficult task which is
connected with the deeper structure of the corresponding Harmonic Analysis
problem. Thus far several Bellman functions have been computed (see
\cite{Burk1}, \cite{Burk2}, \cite{Mel1}, \cite{Sla}, \cite{Sla1}, \cite{Vas},
\cite{Vas1}, \cite{Vas2}). L.Slavin and A.Stokolos \cite{SlSt} linked the
Bellman function computation to solving certain PDE's of the Monge Ampere
type, and in this way they obtained an alternative proof of the Bellman
functions relate to the dyadic maximal operator in \cite{Mel1}. Also in
\cite{Vas2} using the Monge-Ampere equation approach a more general Bellman
function than the one related to the dyadic Carleson imbedding Theorem has be
precisely evaluated thus generalizing the corresponding result in \cite{Mel1}.

However many Bellman functions related to dyadic maximal operators do not obey
the dynamics that make the Monge Ampere approach, or the linearization
approach readily applicable. Such are the cases related to weak $L^{p}$ as
well as more general Lorentz $L^{p,q}$ norms. Recently another approach based
on symmetrization, i.e. decreasing rearrangements, was introduced in
\cite{Mel2} and then refined in \cite{Nik} giving results as the computation
of the Bellman functions related to mixed local $L^{p}\rightarrow L^{q}$
estimates (see \cite{Mel2}) the determination of sharp constants in
$L^{p,\infty}\rightarrow L^{p,\infty}$ and in more general Lorentz
$L^{p,q}\rightarrow L^{p,q}$ norm estimates for dyadic maximal operators (see
\cite{Nik}) and also another proof of the result in \cite{Mel1} (see
\cite{Mel4}). \ This method is based on the following Theorem essentially
proved in \cite{Nik} (see also \cite{Mel2} for a weaker version) and it refers
to the maximal operator $M_{\mathcal{T}}$ defined for any nonatomic
probability space $(X,\mu)$, equipped with any tree-like family $\mathcal{T}$
with (see \cite{Mel1}):

\begin{theorem}
Let $G:[0,+\infty)\rightarrow\lbrack0,+\infty)$ be non-decrasing,
$h:(0,1]\rightarrow$ $\mathbb{R}^{+}$ be any locally integrable function. Then
for any nonatomic probability space $(X,\mu)$, equipped with any tree-like
family $\mathcal{T}$ , for any non-increasing right continuous integrable
function $g:(0,1]\rightarrow$ $\mathbb{R}^{+}$ and any $k\in(0,1]$, the
following equality holds:%
\begin{gather*}
\sup\left\{  \int_{0}^{k}G[(M_{\mathcal{T}}\phi)^{\ast}(t)]h(t)dt:\phi\text{
measurable on }X\text{ with }\phi^{\ast}=g\right\}  =\\
=\int_{0}^{k}G\left(  \frac{1}{t}\int_{0}^{t}g(u)du\right)  h(t)dt.
\end{gather*}

\end{theorem}

Here $\phi^{\ast}$ denotes the equimeasurable decreasing rearrangement of the
measurable function $\phi:X\rightarrow\mathbb{R}$ which is defined on $(0,1]$
since $X$ is a probability space. For completeness we will give here a simpler
proof of the above Theorem. This enables as to reduce the problem of
determining a Bellman type function for the local tree maximal operator
$M_{\mathcal{T}}$ to a problem of a similar nature but on $(0,1]$ and for the
local Hardy operator $\mathcal{H}(g)(t)=\frac{1}{t}\int_{0}^{t}g$ acting on
decreasing functions $g$. This idea applied to convex $G$'s has lead to the
determination of the Bellman functions
\begin{gather}
\mathcal{B}_{p,q}(F,f,L)=\sup{\huge \{}\dfrac{1}{\left\vert Q\right\vert }%
\int_{E}(M_{d}\phi)^{p}:\operatorname{Av}_{Q}(\phi^{p})=F,\operatorname{Av}%
_{Q}(\phi)=f,\nonumber\\
\text{ \ \ \ \ \ \ \ \ \ \ \ \ \ \ \ \ \ }\sup\limits_{R:Q\subseteq
R}\operatorname{Av}_{R}(\phi)=L,E\subseteq Q,\left\vert E\right\vert
=k{\huge \}}\label{i6a}%
\end{gather}
whenever $1\leq q<p$ which are given implicitly via certain solutions of
related ODE's (see \cite{Mel2}). However Theorem 1 (see (\cite{Nik}) allows us
to treat problems of more general nature and the purpose of this paper is to
present certain applications of this method in the case of Lorentz type estimates.

Our first application is related to multiple weak-type estimates and is
described in the following

\begin{theorem}
Given the real numbers $f,F_{1},...,F_{m}>0$ and $p_{1},...,p_{m}>1$ with
$f\leq\min\{\left(  \frac{p_{j}}{p_{j}-1}F_{j}\right)  ^{1/p_{j}}:1\leq j\leq
m\}$ and given any nondecreasing $G:[0,+\infty)\rightarrow\lbrack0,+\infty)$
and $h:(0,1]\rightarrow$ $\mathbb{R}^{+}$ be any locally integrable function,
we have for the following Bellman type function
\begin{align*}
\mathcal{B}_{G,h,p_{1},...,p_{m}}^{\mathcal{T}}(F_{1},...F_{m},f,k) &
=\sup\{\int_{0}^{k}G[(M_{\mathcal{T}}\phi)^{\ast}(t)]h(t)dt:\phi\geq0\text{
measurable on }X\text{ }\\
\text{with }\left\Vert \phi\right\Vert _{1} &  =f,\text{ }\left\Vert
\phi\right\Vert _{p_{1},\infty}\leq F_{1},...,\left\Vert \phi\right\Vert
_{p_{m},\infty}\leq F_{m}\}
\end{align*}
the equality%
\begin{equation}
\mathcal{B}_{G,h,p_{1},...,p_{k}}^{\mathcal{T}}(F_{1},...F_{k},f,k)=\int
_{0}^{k}G\left(  \frac{1}{t}\int_{0}^{\min(t,\sigma\}}\min_{1\leq j\leq
m}\left(  \frac{F_{j}}{u}\right)  ^{1/p_{j}}du\right)  h(t)dt\label{i8}%
\end{equation}
where $\sigma$ is defined by the equality%
\begin{equation}
\int_{0}^{\sigma}\min_{1\leq j\leq m}\left(  \frac{F_{j}}{u}\right)
^{1/p_{j}}du=f.\label{i9}%
\end{equation}

\end{theorem}

Using the above Theorem we find the $L^{p,\infty}\rightarrow L^{q,r}$ Lorentz
type Bellman function for the maximal operator. To make the result more
readable let us denote by $p^{\prime},q^{\prime}$ the dual exponents of
$p,q>1$ (so $p^{\prime}=\frac{p}{p-1}$)

\begin{theorem}
Given $1<q<p$ and $r>0$ the Bellman function:%
\begin{gather}
\mathcal{B}_{(p,\infty),(q,r)}^{\mathcal{T}}(F,f,L)=\sup\{\left\Vert
\max(M_{\mathcal{T}}\phi,L)\right\Vert _{L^{q,r}(X,\mu)}^{r}:\phi\geq0\text{
is measurable with}\nonumber\\
\left\Vert \phi\right\Vert _{L^{1}(X,\mu)}=f,\text{ }\left\Vert \phi
\right\Vert _{L^{p,\infty}(X\mu)}^{p}=F\} \label{i10}%
\end{gather}
defined for $0<f<p^{\prime}F^{1/p}$ and $f\leq L$ is given by%
\begin{equation}
\mathcal{B}_{(p,\infty),(q,r)}^{\mathcal{T}}(F,f,L)=\left\{
\begin{array}
[c]{l}%
\dfrac{q(p-1)q^{\prime}}{r(p-q)}(p^{\prime})^{p^{\prime}r/q^{\prime}}%
f^{\frac{r(p-q)}{q(p-1)}}F^{\frac{r(q-1)}{q(p-1)}}+\dfrac{q}{r}L^{r}-\dfrac
{q}{r}q^{\prime}f^{\frac{r}{q}}L^{\frac{r}{q}(q-1)}\text{ }\\
\text{when }L\leq(p^{\prime})^{p^{\prime}}\left(  \frac{F}{f}\right)
^{1/(p-1)}\\
\dfrac{q(p^{\prime})^{rp/q}}{r(\frac{p}{q}-1)}F^{\frac{r}{q}}L^{r(1-\frac
{p}{q})}+L^{r}\text{ }\\
\text{when }L\geq(p^{\prime})^{p^{\prime}}\left(  \frac{F}{f}\right)
^{1/(p-1)}.
\end{array}
\right.  \label{i11}%
\end{equation}

\end{theorem}

The proofs of the above two Theorems are given in section 3.

Next we define the Bellman function related to a Lorentz $L^{p,q}\rightarrow
L^{p,q}$ type estimate for the (martingale) maximal operator, where $p,\dot
{q}>1$ are arbitrary%
\begin{gather}
\mathcal{BL}_{p,q}^{\mathcal{T}}(F,f)=\sup\{\left\Vert M_{\mathcal{T}}%
\phi\right\Vert _{L^{p,q}(X,\mu)}^{q}:\phi\geq0\text{ is measurable
with}\nonumber\\
\left\Vert \phi\right\Vert _{L^{1}(X,\mu)}=f,\text{ }\left\Vert \phi
\right\Vert _{L^{p,q}(X\mu)}^{q}=F\}.\label{i12}%
\end{gather}
In \cite{Nik} it has been proved that $M_{\mathcal{T}}$ satisfies an
$L^{p,q}\rightarrow L^{p,q}$ estimate with best constant $p^{\prime}$. Here we
will determine the exact form of the corresponding Bellman function
(\ref{i12}). We have.

\begin{theorem}
The Bellman function (\ref{i12}) is defined for all pairs $(F,f)$ with (i)
$0<f^{q}\leq\left(  \frac{p^{\prime}}{q^{\prime}}\right)  ^{q-1}F$ if $1<p\leq
q$ and (ii) $0<f^{q}\leq\frac{q}{p}F$ if $1<q<p$ and in both cases it is given
by%
\begin{equation}
\mathcal{BL}_{p,q}^{\mathcal{T}}(F,f)=\left(  \frac{p^{\prime}}{q^{\prime}%
}\right)  ^{q}\omega_{q}\left(  \left(  \frac{q^{\prime}}{p^{\prime}}\right)
^{q-1}\frac{f^{q}}{F}\right)  ^{q}F\text{.}\label{i13}%
\end{equation}

\end{theorem}

Here $\omega_{q}:$ $[0,1]\rightarrow\lbrack1,q^{\prime}]$ is the inverse
function of $H_{q}(z)=-(q-1)z^{q}+qz^{q-1}$ (defined on $[1,q^{\prime}]$) thus
the same function as the one appearing in the Bellman functions of the usual
$L^{p}$ norms. Note though that in the case $1<q<p$ only a restriction of
$H_{q}$ is inverted (see the proof of this theorem). Also we note that the
case $q=1$ could be inferred from this Theorem but it is easy to see that
since $\int_{0}^{1}t^{\frac{1}{p}-2}\int_{0}^{t}g(u)dudt=p^{\prime}$ $\int
_{0}^{1}(t^{\frac{1}{p}-1}-1)g(t)dt$ that from Theorem 1 $\mathcal{BL}%
_{p,1}^{\mathcal{T}}(F,f)=p^{\prime}(F-f)$. \ In section 4 we will prove
Theorem 4.

\section{Trees and maximal operators}

As in \cite{Mel1} we let $(X,\mu)$ be a nonatomic probability space (i.e.
$\mu(X)=1$). Two measurable subsets $A$, $B$ of $X$ will be called almost
disjoint if $\mu(A\cap B)=0$. Then we give the following.

\begin{definition}
A set $\mathcal{T}$ of measurable subsets of $X$ will be called a tree if the
following conditions are satisfied:

(i) $X\in\mathcal{T}$ \ and for every $I\in\mathcal{T}$ \ we have $\mu(I)>0$.

(ii) For every $I\in\mathcal{T}$ \ there corresponds a finite subset
$\mathcal{C}(I)\subseteq\mathcal{T}$ \ containing at least two elements such that:

\qquad(a) the elements of $\mathcal{C}(I)$ are pairwise almost disjoint
subsets of $I$,

\qquad(b) $I=\bigcup\mathcal{C}(I)$.

(iii) $\mathcal{T}=\bigcup_{m\geq0}\mathcal{T}_{(m)}$ where $\mathcal{T}%
_{(0)}=\{X\}$ and $\mathcal{T}_{(m+1)}=\bigcup_{I\in\mathcal{T}_{(m)}%
}\mathcal{C}(I)$.

(iv) We have $\lim\limits_{m\rightarrow\infty}\sup\limits_{I\in\mathcal{T}%
_{(m)}}\mu(I)=0$.{}
\end{definition}

By removing the measure zero exceptional set $E(\mathcal{T})=\bigcup
_{I\in\mathcal{T}}\bigcup_{\substack{J_{1},J_{2}\in\mathcal{C}(I)\\J_{1}\neq
J_{2}}}(J_{1}\cap J_{2})$ we may replace the almost disjointness above by disjointness.

Now given any tree $\mathcal{T}$ we define the maximal operator associated to
it as follows
\begin{equation}
M_{\mathcal{T}}\phi(x)=\sup\left\{  \frac{1}{\mu(I)}\int_{I}\left\vert
\phi\right\vert d\mu:x\in I\in\mathcal{T}\right\}  \label{t1}%
\end{equation}
for every $\phi\in L^{1}(X,\mu)$.

The above setting can be used not only for the dyadic maximal operator but
also for the maximal operator on martingales, hence many of the results here
can be viewed as generalizations and refinements of the classical Doob's inequality.

The following Lemma has been proved in \cite{Mel1} and provides the basis of
constructing examples that show sharpness.

\begin{lemma}
For every $I\in\mathcal{T}$ and every $\alpha$ such that $0<\alpha<1$ there
exists a subfamily $\mathcal{F}(I)\subseteq\mathcal{T}$ \ consisting of
pairwise almost disjoint subsets of $I$ such that
\begin{equation}
\mu(\bigcup_{J\in\mathcal{F}(I)}J)=\sum_{J\in\mathcal{F}(I)}\mu(J)=(1-\alpha
)\mu(I)\text{.} \label{t2}%
\end{equation}

\end{lemma}

Then we have the following Lemma which give the one side of Theorem 1.

\begin{lemma}
\bigskip Given any nonnegative integrable $\phi$ on $X$ we have$^{\ast}$
\begin{equation}
(M_{\mathcal{T}}\phi)^{\ast}(t)\leq\frac{1}{t}\int_{0}^{t}\phi^{\ast}(u)du
\label{t3}%
\end{equation}
for every $t\in(0,1)$ and therefore whenever $G,h,k$ are as in Theorem 1
\begin{equation}
\int_{0}^{k}G[(M_{\mathcal{T}}\phi)^{\ast}(t)]h(t)dt\leq\int_{0}^{k}G\left(
\frac{1}{t}\int_{0}^{t}\phi^{\ast}(u)du\right)  h(t)dt. \label{t4}%
\end{equation}

\begin{proof}
Fixing $t\in(0,1)$ let $\alpha=(M_{\mathcal{T}}\phi)^{\ast}(t)=\inf
\{\lambda:\mu(\{M_{\mathcal{T}}\phi\geq\lambda\})\leq t\}$. Then given any
$\lambda<\alpha$ we have $\mu(\{M_{\mathcal{T}}\phi\geq\lambda\})>t$ and using
the decomposition of $\{M_{\mathcal{T}}\phi\geq\lambda\}$ as a disjoint union
of elements $I$ of $\mathcal{T}$ maximal under the condition $\int_{I}\phi
d\mu\geq\lambda\mu(I)$, we conclude that there exists measurable $A\subseteq
X$ with $\mu(A)>t$ and $\int_{A}\phi d\mu\geq\lambda\mu(A)$. But now since
$\phi^{\ast}$ is decreasing we have
\begin{equation}
\frac{1}{t}\int_{0}^{t}\phi^{\ast}\geq\frac{1}{\mu(A)}\int_{0}^{\mu(A)}%
\phi^{\ast}\geq\frac{1}{\mu(A)}\int_{A}\phi d\mu\geq\lambda. \label{t5}%
\end{equation}
This holding for any $\lambda<\alpha$ implies (\ref{t3}).
\end{proof}
\end{lemma}

\bigskip The construction in the next Lemma appears also in \cite{Nik} and
provides the other half of Theorem 1. We include a simpler proof for completeness.

\begin{lemma}
For $G,h,k$ and $g$ as in Theorem 1, there exists a sequence of measurable
functions $\psi_{N}:X\rightarrow$ $\mathbb{R}^{+}$ such that $\psi_{N}^{\ast
}=g$ and
\begin{equation}
\lim\sup\nolimits_{N\rightarrow\infty}\int_{0}^{k}G[(M_{\mathcal{T}}\psi
_{N})^{\ast}(t)]h(t)dt\geq\int_{0}^{k}G\left(  \frac{1}{t}\int_{0}%
^{t}g\right)  h(t)dt\text{.} \label{t6}%
\end{equation}

\begin{proof}
Fixing $\alpha$ with $0<\alpha<1$ and using Lemma 1, we choose for every
$I\in\mathcal{T}$ a family $\mathcal{F}(I)\subseteq\mathcal{T}$ of pairwise
almost disjoint subsets of $I$ such that
\begin{equation}
\sum_{J\in\mathcal{F}(I)}\mu(J)=(1-\alpha)\mu(I)\text{.}\label{e13}%
\end{equation}
Then we define $\mathcal{S}=\mathcal{S}_{\alpha}$ to be the smallest subset of
$\mathcal{T}$ such that $X\in\mathcal{S}$ and for every $I\in\mathcal{S}$,
$\mathcal{F}(I)\subseteq\mathcal{S}$. Next for every $I\in\mathcal{S}$ we
define the set
\begin{equation}
A_{I}=I~\backslash%
{\displaystyle\bigcup\limits_{J\in\mathcal{F}(I)}}
J\label{e14}%
\end{equation}
and note that $\mu(A_{I})=\alpha\mu(I)$ and $I=%
{\displaystyle\bigcup\limits_{_{\substack{J\in\mathcal{S}\\J\subseteq I}}}}
A_{J}$ for every $I\in\mathcal{S}$. Also since $\mathcal{S}=\bigcup_{m\geq
0}\mathcal{S}_{(m)}$ where $\mathcal{S}_{(0)}=\{X\}$ and $\mathcal{S}%
_{(m+1)}=\bigcup_{I\in\mathcal{S}_{(m)}}\mathcal{F}(I)$, we can define
rank$(I)=r(I)$ for $I\in\mathcal{S}$ to be the unique integer $m$ such that
$I\in\mathcal{S}_{(m)}$ and remark that $\sum\limits_{\substack{\mathcal{S}\ni
J\subseteq I\\r(J)=r(I)+m}}\mu(J)=(1-\alpha)^{m}\mu(I)$ for every
$I\in\mathcal{S}$. For any $m\geq0$ let
\begin{equation}
\gamma_{m}=\frac{1}{\alpha(1-\alpha)^{m}}\int_{(1-\alpha)^{m+1}}%
^{(1-\alpha)^{m}}g(u)du\label{e16}%
\end{equation}
and for any $I\in\mathcal{S}_{(m)}$ i.e. rank$(I)=m$, since $\mu$ is nonatomic
we can choose a random variable $R_{I}:A_{I}\rightarrow\lbrack0,+\infty)$ on
the probability space $(A_{I},\frac{1}{\mu(A_{I})}\mu)$ having the same
distribution as the restriction of $g$ on the probability space $((1-a)^{m+1}%
,(1-\alpha)^{m}]$ with measure $\frac{1}{\alpha(1-\alpha)^{m}}d\lambda$
($\lambda$ being Lebesgue measure). Then define%
\[
\phi_{\alpha}(x)=R_{I}(x)\text{ when }x\in A_{I},~I\in\mathcal{S}\text{.}%
\]
For and any $s>0$ the disjointness of the $A_{I}$'s implies that%
\begin{align*}
\mu(\{\phi_{\alpha} &  \geq s\})=\sum_{I\in\mathcal{S}}\mu(\{x\in A_{I}%
:R_{I}(x)\geq s\})=\\
&  =\sum_{m\geq0}\sum_{I\in\mathcal{S}_{(m)}}\frac{\mu(A_{I})}{\alpha
(1-\alpha)^{m}}\left\vert \{t\in((1-a)^{m+1},(1-\alpha)^{m}]:g(t)\geq
s\}\right\vert =\\
&  =\left\vert \{t\in(0,1]:g(t)\geq s\}\right\vert
\end{align*}
hence $\phi_{\alpha}$ and $g$ have the same distribution and since $g$ is
nonincreasing and right continuous on $(0,\mu(X)]$ we conclude that
$\phi_{\alpha}^{\ast}=g$. Moreover
\[
\frac{1}{\mu(A_{I})}\int_{A_{I}}\phi_{a}d\mu=\gamma_{m}%
\]
for every $I\in\mathcal{S}$ with rank$(I)=m$ and thus%
\begin{align*}
\operatorname{Av}_{I}(\phi_{\alpha}) &  =\dfrac{1}{\mu(I)}\int_{I}\phi_{a}%
d\mu=\dfrac{1}{\mu(I)}\sum_{_{\substack{J\in\mathcal{S}\\J\subseteq I}}}%
\int_{A_{J}}\phi_{a}d\mu=\\
&  =\dfrac{\alpha}{\mu(I)}\sum_{\ell\geq0}\gamma_{\ell+\text{rank}(I)}%
\sum_{\substack{\mathcal{S}\ni J\subseteq I\\\text{rank}(J)=\text{rank}%
(I)+\ell}}\mu(J)=\\
&  =\alpha\sum_{\ell\geq0}\gamma_{\ell+m}(1-\alpha)^{\ell}=\frac{1}%
{(1-\alpha)^{m}}\int_{0}^{(1-\alpha)^{m}}g(u)du
\end{align*}
implying that $M_{\mathcal{T}}\phi\geq\frac{1}{(1-\alpha)^{m}}\int
_{0}^{(1-\alpha)^{m}}g(u)du$ on the set $%
{\displaystyle\bigcup\limits_{I\in\mathcal{S}:r(I)=m}}
I$ which has measure $(1-\alpha)^{m}$ and thus $(M_{\mathcal{T}}\phi_{\alpha
})^{\ast}(t)\geq\frac{1}{(1-\alpha)^{m}}\int_{0}^{(1-\alpha)^{m}}g(u)du$ for
every $t\in((1-a)^{m+1},(1-\alpha)^{m})$.

Now with $N$ large, taking $\alpha_{N}=1-(1-k)^{1/N}$ and $\psi_{N}%
=\phi_{\alpha_{N}}$ we have
\[
\int_{0}^{k}G[(M_{\mathcal{T}}\psi_{N})^{\ast}(t)]h(t)dt\geq\sum_{j\geq
0}G\left(  \frac{1}{(1-k)^{1+\frac{m}{N}}}\int_{0}^{(1-k)^{1+\frac{m}{N}}%
}g(u)du\right)  \int_{(1-k)^{1+\frac{m+1}{N}}}^{(1-k)^{1+\frac{m}{N}}}h(t)dt
\]
the last sum converging to $\int_{0}^{k}G\left(  \frac{1}{t}\int_{0}%
^{t}g\right)  h(t)dt$ as $N\rightarrow\infty$ by monotone convergence. This
completes the proof.
\end{proof}
\end{lemma}

\section{The case of weak type conditions}

Here we will prove Theorems 2 and 3. Theorem 2 follows from the following more
general Proposition by taking $R(t)$ to be the decreasing function
$\min_{1\leq j\leq m}\left(  \frac{F_{j}}{t}\right)  ^{1/p_{j}}$.

\begin{proposition}
\bigskip Let $R:(0,1]\rightarrow(0,+\infty)$ be a decreasing, continuous and
integrable function and for any $f$ with $0<f\leq\int_{0}^{1}R(t)dt$ let
$\sigma=\sigma(f)$ be the unique number in $(0,1]$ with $\int_{0}^{\sigma
}R(t)dt=f$. Then for any $(X,\mu,\mathcal{T}),G,h,k$ as in Theorem 1 we have%
\begin{gather}
\sup\left\{  \int_{0}^{k}G[(M_{\mathcal{T}}\phi)^{\ast}(t)]h(t)dt:\phi
\geq0\text{ with }\int_{X}\phi d\mu=f\text{ and  }\phi^{\ast}\leq R\right\}
=\nonumber\\
=\int_{0}^{k}G\left(  \frac{1}{t}\int_{0}^{t}R(u)\chi_{\lbrack0,\sigma
(f))}(u)du\right)  h(t)dt.\label{w1}%
\end{gather}

\begin{proof}
By Theorem the above type supremum but fixing $\phi^{\ast}=g$ is equal to
$\int_{0}^{k}G\left(  \frac{1}{t}\int_{0}^{t}g(u)du\right)  h(t)dt$. But we
must have $g(t)\leq R(t)$ for any $t$ and $\int_{0}^{1}g=\int_{X}\phi
d\mu=f=\int_{0}^{\sigma}R$ we conclude that one the one hand $\int_{0}^{t}g$
$\leq\int_{0}^{t}R$ when $0\leq t\leq\sigma$ and on the other hand $\int
_{0}^{t}g\leq\int_{0}^{1}g=\int_{0}^{\sigma}R=\int_{0}^{t}R(u)\chi
_{\lbrack0,\sigma(f))}(u)du$ when $\sigma<t\leq1$ we get $G\left(  \frac{1}%
{t}\int_{0}^{t}g\right)  \leq G\left(  \frac{1}{t}\int_{0}^{t}R\chi
_{\lbrack0,\sigma)}\right)  $ for all $t$. Thus using the converse implication
in Theorem 1 for the decreasing right continuous function $R\chi
_{\lbrack0,\sigma)}$ completes the proof of (\ref{w1}).
\end{proof}
\end{proposition}

Now to prove Theorem 3 we remark that using Theorem 2 with $G(x)=\max
(x,L)^{r}$, $h(t)=t^{\frac{r}{q}-1}$ and $m=1,p_{1}=p$ that the expression
$\mathcal{B}_{(p,\infty),(q,r)}^{\mathcal{T}}(F,f,L)$ in (\ref{i10}) is equal
to the following expression (actually we get the supremum under $\left\Vert
\phi\right\Vert _{L^{p,\infty}(X\mu)}^{p}\leq F$ but it is easy to see that at
the extremum above we have the equality $\left\Vert \phi\right\Vert
_{L^{p,\infty}(X\mu)}^{p}=F$)%
\begin{equation}
\int_{0}^{1}t^{\frac{r}{q}-1}\max\left(  \frac{1}{t}\int_{0}^{\min(t,\sigma
\}}\left(  \frac{F}{u}\right)  ^{1/p}du,L\right)  ^{r}dt\label{w2}%
\end{equation}
where $\sigma=\sigma_{p}(F,f)$ is given by $\int_{0}^{\sigma}\left(  \frac
{F}{u}\right)  ^{1/p}du=f$ thus%
\begin{equation}
\sigma=\left(  \frac{f}{p^{\prime}F^{1/p}}\right)  ^{p^{\prime}}.\label{w3}%
\end{equation}
Next note that using (\ref{w3})
\begin{equation}
\frac{1}{t}\int_{0}^{\min(t,\sigma\}}\left(  \frac{F}{u}\right)
^{1/p}du=p^{\prime}F^{1/p}\frac{1}{t}\min(t,\sigma)^{1-\frac{1}{p}}%
=\min\left(  p^{\prime}\left(  \frac{F}{t}\right)  ^{1/p},\frac{f}{t}\right)
\label{w4}%
\end{equation}
and so we have to compute the integral $\int_{0}^{1}t^{\frac{r}{q}-1}%
\Sigma(t)^{r}dt$ where $\Sigma(t)$ is given by%
\begin{equation}
\Sigma(t)=\max\left(  \min\left(  p^{\prime}\left(  \frac{F}{t}\right)
^{1/p},\frac{f}{t}\right)  ,L\right)  .\label{w5}%
\end{equation}
Observing that $L>\min\left(  p^{\prime}\left(  \frac{F}{t}\right)
^{1/p},\frac{f}{t}\right)  $ if and only if $t>\frac{f}{L}$ or $t>\left(
\frac{p^{\prime}}{L}\right)  ^{p}F$ we consider two cases:

\textbf{Case 1} If $\frac{f}{L}\leq\left(  \frac{p^{\prime}}{L}\right)  ^{p}F$
that is $L\leq(p^{\prime})^{p^{\prime}}\left(  \frac{F}{f}\right)  ^{\frac
{1}{p-1}}$ we easily get $\frac{f}{L}\geq\left(  \frac{f}{p^{\prime}F^{1/p}%
}\right)  ^{p^{\prime}}=\sigma$ and therefore we have%
\[
\Sigma(t)=\left\{
\begin{array}
[c]{l}%
p^{\prime}\left(  \frac{F}{t}\right)  ^{1/p}\text{ when }0\leq t\leq\sigma\\
\frac{f}{t}\text{ when }\sigma<t\leq\frac{f}{L}\\
L\text{ when }\frac{f}{L}<t\leq1
\end{array}
\right.
\]
and then computing the corresponding integral $\int_{0}^{1}t^{\frac{r}{q}%
-1}\Sigma(t)^{r}dt$ we get the upper half in (\ref{i12}).

\textbf{Case 2} If $\frac{f}{L}>\left(  \frac{p^{\prime}}{L}\right)  ^{p}F$
that is $L>(p^{\prime})^{p^{\prime}}\left(  \frac{F}{f}\right)  ^{\frac
{1}{p-1}}$ we easily get $\left(  \frac{p^{\prime}}{L}\right)  ^{p}F<\left(
\frac{f}{p^{\prime}F^{1/p}}\right)  ^{p^{\prime}}=\sigma$ and therefore we
have%
\[
\Sigma(t)=\left\{
\begin{array}
[c]{l}%
p^{\prime}\left(  \frac{F}{t}\right)  ^{1/p}\text{ when }0\leq t\leq\left(
\frac{p^{\prime}}{L}\right)  ^{p}F\\
L\text{ when }\left(  \frac{p^{\prime}}{L}\right)  ^{p}F<t\leq1
\end{array}
\right.
\]
and then computing the corresponding integral $\int_{0}^{1}t^{\frac{r}{q}%
-1}\Sigma(t)^{r}dt$ we get the lower half in (\ref{i12}).

These cases complete the proof of Theorem 3.

\bigskip

\section{\bigskip Proof of Theorem 4}

In view of Theorem 1 and by setting $\phi^{\ast}=g$ it suffices to determine
the supremum of the expression $\Delta(g)=\int_{0}^{1}(t^{\frac{1}{p}-1}%
\int_{0}^{t}g(u)du)^{q}\frac{dt}{t}$ when $g$ runs over all nonnegative
decreasing right continuous functions on $(0,1]$ satisfying $\int_{0}%
^{1}g(t)dt=f$ and $\int_{0}^{1}(t^{\frac{1}{p}}g(t))^{q}\frac{dt}{t}=F$.
Considering first any bounded such function $g$ we compute by integration by
parts%
\begin{gather*}
\int_{0}^{1}t^{q(\frac{1}{p}-1)}(\int_{0}^{t}g(u)du)^{q-1}g(t)dt=\frac{1}%
{q}\int_{0}^{1}t^{q(\frac{1}{p}-1)}[(\int_{0}^{t}g(u)du)^{q}]^{\prime}dt=\\
=\frac{1}{q}(\int_{0}^{1}g(u)du)^{q}+(1-\frac{1}{p})\int_{0}^{1}(t^{\frac
{1}{p}-1}\int_{0}^{t}g(u)du)^{q}\frac{dt}{t}=\frac{f^{q}}{q}+\frac
{1}{p^{\prime}}\Delta(g)\text{.}%
\end{gather*}
But using Young's inequality $xy\leq\frac{x^{q}}{q}+\frac{y^{q^{\prime}}%
}{q^{\prime}}$ in the first integral as follows,where $\gamma=\frac{1}%
{p}-\frac{1}{q}$ and $\lambda>0$ will be determined later%
\begin{gather*}
\int_{0}^{1}t^{q(\frac{1}{p}-1)}(\int_{0}^{t}g(u)du)^{q-1}g(t)dt=\int_{0}%
^{1}t^{-\gamma+q(\frac{1}{p}-1)}(\frac{1}{\lambda^{1/(q-1)}}\int_{0}%
^{t}g(u)du)^{q-1}(\lambda g(t)t^{\gamma})dt\leq\\
\leq\frac{1}{q}\int_{0}^{1}\lambda^{q}g(t)^{q}t^{\gamma q}dt+\frac
{1}{q^{\prime}}\int_{0}^{1}\lambda^{-q^{\prime}}t^{[-\gamma+q(\frac{1}%
{p}-1)]q^{\prime}}(\int_{0}^{t}g(u)du)^{q}dt=\\
=\frac{\lambda^{q}}{q}\int_{0}^{1}g(t)^{q}t^{\frac{q}{p}-1}dt+\frac
{\lambda^{-q^{\prime}}}{q^{\prime}}\int_{0}^{1}t^{\frac{q}{p}-q-1}(\int
_{0}^{t}g(u)du)^{q}dt=\frac{\lambda^{q}}{q}F+\frac{\lambda^{-q^{\prime}}%
}{q^{\prime}}\Delta(g)\text{.}%
\end{gather*}
Therefore we have by writing $\delta=\frac{p^{\prime}}{q^{\prime}}$ and taking
$\lambda^{q^{\prime}}=(\beta+1)\delta,~\beta>0$ and using the above
inequalities that
\[
\frac{\beta}{p^{\prime}(\beta+1)}\Delta(g)=\left(  \frac{1}{p^{\prime}}%
-\frac{\lambda^{-q^{\prime}}}{q^{\prime}}\right)  \Delta(g)\leq\frac
{\lambda^{q}F-f}{q}%
\]
and so%
\begin{equation}
\Delta(g)\leq\delta\frac{\beta+1}{\beta}\frac{(\beta+1)^{q-1}\delta
^{q-1}F-f^{q}}{q-1}.\label{l1}%
\end{equation}
Next, given an arbitrary $g$, the above estimate can be used for the
truncations $g_{M}=\min(g,M)$ and $F,f$ replaced by the corresponding
quantities for $g_{M}$ and then take $M\rightarrow+\infty$ and use monotone
convergence to infer that (\ref{l1}) holds for the general nonnegative
decreasing right continuous function on $(0,1]$ satisfying $\int_{0}%
^{1}g(t)dt=f$ and $\int_{0}^{1}(t^{\frac{1}{p}}g(t))^{q}\frac{dt}{t}=F$.

As has been also remarked in \cite{Mel1} it is easy to see that the right hand
side of (\ref{l1}) is minimized when $\beta$ satisfies the equation
$H_{q}(\beta+1)=\frac{f^{q}}{\delta^{q-1}F}\leq1$ (which is well known that it
is $\leq1$, but also follows from (\ref{l1}) by taking $\beta\rightarrow0^{+}%
$) and then for this value of $\beta$ the right hand side of (\ref{l1})
becomes $\delta^{q}\omega_{q}\left(  \frac{f^{q}}{\delta^{q-1}F}\right)
^{q}F$. This proves the inequality $\mathcal{BL}_{p,q}^{\mathcal{T}}(F,f)\leq$
$\delta^{q}\omega_{q}\left(  \frac{f^{q}}{\delta^{q-1}F}\right)  ^{q}F$.

Now we consider the continuous positive decreasing function%
\begin{equation}
g_{\alpha}(t)=f(1-\alpha)t^{-\alpha}\label{l2}%
\end{equation}
where $0\leq\alpha<1$. Clearly $\int_{0}^{1}g_{\alpha}(t)dt=f$ and since
$\frac{1}{t}\int_{0}^{t}g_{\alpha}(u)du=\frac{g_{\alpha}(t)}{1-\alpha}$ for
all $t\in(0,1]$ we have $\Delta(g_{\alpha})=\left(  \frac{1}{1-\alpha}\right)
^{q}\int_{0}^{1}(t^{\frac{1}{p}-1}g_{\alpha}(t))^{q}\frac{dt}{t}$. The
condition $\int_{0}^{1}(t^{\frac{1}{p}}g_{\alpha}(t))^{q}\frac{dt}{t}=F$ is
then equivalent to the equation in $\alpha$%
\begin{equation}
\frac{p}{q}\frac{(1-\alpha)^{q}}{1-\alpha p}=\frac{F}{f^{q}}\text{.}\label{l3}%
\end{equation}
Consider the function $w(\alpha)=\frac{p}{q}\frac{(1-\alpha)^{q}}{1-\alpha p}$
defined on $0\leq\alpha<\frac{1}{p}$. \ We have $w^{\prime}(\alpha
)=\frac{p(q-1)\alpha+p-q}{1-\alpha}w(\alpha)$, $w(0)=\frac{p}{q}$ and
$\lim\limits_{\alpha\rightarrow(1/p)^{+}}w(\alpha)=+\infty$. Now consider the
following cases:

\textbf{Case 1 }$1<p\leq q.$In this case the function $w$ has a minimum at
$\alpha_{0}=\frac{q-p}{p(q-1)}<\frac{1}{p}$ and is strictly increasing on
$[\alpha_{0},\frac{1}{p})$ and since also $\alpha_{0}=1-\frac{q^{\prime}%
}{p^{\prime}}=1-\frac{1}{\delta}$ which gives $w(\alpha_{0})=\delta^{q-1}$,
its range is $[\delta^{q-1},+\infty)$. This implies that the domain of
$\mathcal{BL}_{p,q}^{\mathcal{T}}(F,f)$ in this case consists of all $F,f$
with $0<f^{q}\leq\delta^{q-1}F$ as asserted in Theorem 4, and with such a pair
$(F,f)$ there exists a unique $\alpha=\alpha(F,f)$ in the interval
$[\alpha_{0},\frac{1}{p})$ such that $w(a)=\frac{F}{f^{q}}$. Then for this
$\alpha$ we have $\int_{0}^{1}g_{\alpha}(t)dt=f$, $\int_{0}^{1}(t^{\frac{1}%
{p}}g_{\alpha}(t))^{q}\frac{dt}{t}=F$ and $\Delta(g_{\alpha})=\left(  \frac
{1}{1-\alpha}\right)  ^{q}F$ and since $\frac{1}{1-\alpha}\geq\frac
{1}{1-\alpha_{0}}=\delta$ we may write $\frac{1}{1-\alpha}=\delta z$ with
$z\geq1$ (and also $z<\frac{1}{\delta}p^{\prime}=q^{\prime}$) and then it is
easy to see that (\ref{l3}) transforms into $H_{q}(z)=\frac{f^{q}}%
{\delta^{q-1}F}$ and so $z=\omega_{q}\left(  \frac{f^{q}}{\delta^{q-1}%
F}\right)  $ giving that $\Delta(g_{\alpha})=\delta^{q}\omega_{q}\left(
\frac{f^{q}}{\delta^{q-1}F}\right)  ^{q}F$ and thus proving Theorem 4 in this case.

\textbf{Case 2 }$1<q<\dot{p}.$In this case the function $w$ has positive
derivative hence it is one to one and its range is $[\frac{p}{q},+\infty)$. On
the other hand for any nonnegative decreasing right continuous function $g$ on
$(0,1]$ satisfying $\int_{0}^{1}g(t)dt=f$ and $\int_{0}^{1}(t^{\frac{1}{p}%
}g(t))^{q}\frac{dt}{t}=F$, Chebyshev's inequality (applicable since
$t^{\frac{q}{p}-1}$ is here decreasing) gives%
\[
F=\int_{0}^{1}g(t)^{q}t^{\frac{q}{p}-1}dt\geq\int_{0}^{1}g(t)^{q}dt\int
_{0}^{1}t^{\frac{q}{p}-1}dt\geq(\int_{0}^{1}g(t)dt)^{q}\int_{0}^{1}t^{\frac
{q}{p}-1}dt=f^{q}\frac{p}{q}%
\]
therefore on the one hand this proves that the domain of $\mathcal{BL}%
_{p,q}^{\mathcal{T}}(F,f)$ in this case consists of all $F,f$ with
$0<f^{q}\leq\frac{q}{p}F$ as asserted in Theorem 4, and on the other hand
given any pair $(F,f)$ satisfying $\frac{F}{f^{q}}\geq\frac{p}{q}$ in
(\ref{l3}) there exists a unique $\alpha=\alpha(F,f)$ in the interval
$[0,\frac{1}{p})$ such that $w(a)=\frac{F}{f^{q}}$. Since in this case
$\delta=\frac{p^{\prime}}{q^{\prime}}<1$ we may write $\frac{1}{1-\alpha
}=\delta z$ with $z\geq\frac{1}{\delta}>1$ (and also $z<\frac{1}{\delta
}p^{\prime}=q^{\prime}$) and then as in case 1 we get for this $\alpha$,
$\int_{0}^{1}g_{\alpha}(t)dt=f$, $\int_{0}^{1}(t^{\frac{1}{p}}g_{\alpha
}(t))^{q}\frac{dt}{t}=F$ and $\Delta(g_{\alpha})=\delta^{q}\omega_{q}\left(
\frac{f^{q}}{\delta^{q-1}F}\right)  ^{q}F$. The only difference here is that
only the restriction of $H_{q}$ on $[\frac{q^{\prime}}{p^{\prime}},q^{\prime
}]$ is inverted. These complete the proof of Theorem 4.

\end{document}